\documentclass[10pt,a4paper]{article}

\usepackage{amsmath, amsfonts, amssymb}

    \long\def\symbolfootnote[#1]#2{\begingroup%
    \def\thefootnote{\fnsymbol{footnote}}\footnote[#1]{#2}\endgroup} 

\newtheorem{thm}{Theorem}
\newtheorem{lem}[thm]{Lemma}
\newtheorem{prop}[thm]{Proposition}

\newcommand{\A}{A}

\newcommand{\CC}{{\mathbb C}}
\newcommand{\eps}{\varepsilon}
\newcommand{\HH}{{\cal H}}
\newcommand{\II}{{\mathcal I}}
\newcommand{\III}{{\mathcal J}}

\newcommand{\MM}{{\mathbb M}}
\newcommand{\OO}{{\mathcal O}}
\newcommand{\V}{{\mathcal V}}
\newcommand{\rank} {\, {\textrm{rank}} \, }

\newcommand{\EE}{{\mathcal E}}
\newcommand{\D}{{\mathcal D}}
\newcommand{\rr}{{\mathbb R}}
\newcommand{\N}{{\mathbb N}}
\newcommand{\NN}{{\mathcal N}}

\begin{document}
\rm
%%%%%%%%%%%%%%%%%%%%%%%%%%%%%%%%%%%%%%%%%%%%
%%%%%%%%%%%%%%%%%%%%%%%%%%%%%%%%%%%%%%%%%%%%
%%%%%%%%%%%%%%%%%%%%%%%%%%%%%%%%%%%%%%%%%%%%

%%%%%%%%%% Beginning of the text

\centerline{\Large\bf Bounds for Entropy Numbers} 
\centerline{ \Large\bf for Some Critical Operators}
\bigskip

\centerline{M.A. Lifshits\footnote{The work supported by the RFBR-DFG grant 09-01-91331 
%%"Geometry and asymptotics of random structures"
and by RFBR grant 09-01-12180-ofi\_m. 
%%"Development of fundamental methods
%%of analysis and approximation of probabilistic distributions with applications"
}
}

\centerline{\it St.Petersburg State University}

\bigskip

\begin{abstract}
We provide upper bounds for entropy numbers for two types of operators:
summation operators on binary trees and integral operators of Volterra type. 
Our efforts are concentrated on the critical cases where none of known methods works.
Therefore, we develop a method which seems to be completely new and probably merits 
further applications.
\symbolfootnote[0]{\noindent{\bf AMS Subject Classification:}\ 47B06.}
\symbolfootnote[0]{\noindent{\bf Key words and phrases:}\ entropy numbers, integral operators, operators on trees.}
\end{abstract}
\medskip

\section*{Introduction}
 \label{s:intro1}
We will investigate the entropy numbers of certain linear operators. Recall that \ for a set $A$ in a metric space,
its\ $n$-th diadic entropy number  is the infimum of $\eps>0$ such that $A$ admits a covering by $2^{n-1}$
balls of radius $\eps$. Moreover, given a compact linear operator $V:X\to Y$ acting from one normed space to
another, its entropy numbers $e_n(V)$ are defined as those of V-image of the unit ball of $X$.
The entropy numbers along with other measures of compactness such as approximation numbers and Kolmogorov
numbers play extremely important role in operator theory and its applications. We refer to classical
monographs \cite{CS} and \cite{ET} for further details and references.  

This work originates from a question of M. Lacey and W. Linde. They investigated entropy numbers for linear 
Volterra operators with relatively bad compactness properties and discovered that two types of the behavior
of entropy numbers are possible \cite{La}, \cite{Lin} (see more details in Section \ref{s:integral}). 
On a certain boundary separating the two cases their methods did not apply
and the problem remained open. Further hard efforts convinced us that the remaining case can not be settled 
by a rich variety of traditional methods.  Therefore a new technique is required. It turned out that this new technique
could be cleanly elaborated and better explained if we replace  Volterra operator 
by analogous summation operator on the binary tree. This class of operators is quite simple and natural but
it is absolutely not investigated (its properties will be a subject of a separate work). 
Therefore, we start with consideration of summation operators and first prove our estimate in this case.
Notice that the trees appear naturally in the study of functional spaces because Haar base and other
similar wavelet bases indeed have a structure close to that of a binary tree.

In the last section, we reproduce the same approach for the integral operator considered by Lacey and Linde.

Here is our main point: in most part of classical methods to evaluate entropy numbers $e_n(\V)$ of an operator
$\V$ one approximates $\V$ with a finite rank operator depending on $n$. Contrary to this, we approximate $\V$ with a
{\it family} of finite rank operators indexed by some finite set of\ "essential trees", a notion introduced
in this article.

\section{Introduction to tree summation operators}
\label{s:intro2}
We consider a tree $T$ and its levels $\{T_l\}$, 
$l=0,1,\dots $ such that the level $T_0$ consists of the single node
(the tree root) and the level $T_{l+1}$ is the set of all direct offsprings of nodes that belong to
$T_l$.   

We denote $\OO_T(t)=\OO(t)$ the set of all direct and indirect offsprings 
of a node $t\in T$ including $t$ itself  and let $\OO_l(t)=\OO(t)\cap T_l$.
If $t\in T_l$, we write $|t|=l$.
If $u\in \OO(t)$ we write $u\succeq t$ and $t\preceq u$. 
The strict inequalities have the same meaning with additional assumption $u\not=t$.

For any element $\mu\in \ell_1(T)$ and any $t\in T$ we denote the mass and variation
at $t$ as
\[
s_\mu(t)=\sum_{u\succeq t} \mu(u),\qquad
||\mu||(t)=\sum_{u\succeq t} |\mu(u)|.
\]
Clearly, for any $t\in T$
\begin{equation} \label{s1}
   s_\mu(t)\le ||\mu||(t)
\end{equation}
 and for any $t\in T_l$ and any $m\ge l$ we have 
\begin{equation} \label{s2}
    \sum_{u\in \OO_m(t)} ||\mu||(u) \le ||\mu||(t).
\end{equation}

Now  assume that $T$ is equipped with a non-negative weight
$W=\{w(t)\}_{t\in T}$.

The weight $W$ gives rise to the following simple 
weighted summation operator,
$\tilde \V: \ell_2(T)\to \ell_\infty(T)$  given by  
\[
(\tilde \V\mu) (t)= \sum_{u\preceq t} \sqrt{w(u)} f(u), \qquad t\in T.
\]
where the summation is actually taken over the branch leading from the root to the node $t$.
By technical reasons, we will investigate a slightly different form of
this operator. Namely, let us introduce a pair of 
dual {\it tree-summation operators},
$\V:\ell_2(T,W)\to \ell_\infty(T)$ and $\V^*: \ell_1(T)\to \ell_2(T,W)$ defined by
\begin{equation} \label{v_op}
(\V f) (t)= \sum_{u\preceq t} w(u) f(u), \qquad t\in T,
\end{equation}
and
\begin{equation} \label{vstar_op}
(\V^*\mu) (t)=s_\mu(t)= \mu(\OO(t)), \qquad t\in T,
\end{equation}
respectively. It is easy to see that 
\[
||\V||^2= ||\V^*||^2 = \sup_{t\in T} \ \sum_{u\preceq t} w(u).
\]
It is also clear that the operators $\V$ and $\tilde \V$ are isomorphic. 
We have chosen the representation (\ref{v_op}) because of the simple
form of the operator $\V^*$, the one we will really handle.

\section{The entropy of a summation operator on the binary tree}
\label{s:binary}

In this section  we consider a {\it binary} tree $T$ with levels $\{T_l\}$, 
$l=0,1,\dots $ such that the level $T_0$ consists of the single node
(the tree root) and every node of level $T_l$ generates $2$ offsprings
in $T_{l+1}$. Note that $|T_l| ={2^l}$.

The  weight  $W=\{w(t)\}_{t\in T}$ is defined by
\begin{equation}\label{bw}
     w(t)=(1+|t|)^{-\beta}, \qquad t\in T,\ \ \beta>1.
\end{equation}

\subsection{Regular case}

\begin{thm} \label{t2noncrit} Let $\beta>1$ and let the weight $W$ be given by $(\ref{bw})$.
Consider the linear operator
$\V^*:\ell_1(T)\to \ell_2(T,W)$
defined by $(\ref{vstar_op})$. There exist numeric constants $C_1, C_2$ depending on
$\beta$ such that
for all positive integers $n$ we have the following  bounds for its entropy
numbers
\[
   {C_1}{n^{-\frac{\beta-1}2}} \ \le \ e_n(\V^*) \ \le \ {C_2}{n^{-\frac{\beta-1}2}}\ , \qquad 1<\beta<2 ; 
\]
\[
   {C_1(\ln n)^{1-\beta/2}}{n^{-1/2}} \  \le \  e_n(\V^*) \ \le \ {C_2(\ln n)^{1-\beta/2}}{n^{-1/2}}\ , \qquad \beta>2 ; 
\]
\[
   {C_1}{n^{-1/2}} \ \le \ e_n(\V^*) \ \le \ {C_2 \ (\ln n) }{n^{-1/2}}\ , \qquad \beta=2 . 
\]
\end{thm} 

{\bf Proof.} \
{\it Upper bound}. Consider the set $D=\{\V^*\delta_t, t\in T\}\subset \ell_{2,W}(T)$, where as usual
$\delta_t$ denotes the delta-function at point $t$, i.e $\delta_t(u)=1_{\{u=t\}}$.
Recall that 
\[
   (\V^*\delta_t)(u)=1_{\{u  \preceq t \}}.
\]
It is easy to establish an upper bound for diadic entropy of $D$. Indeed, take a net
$D_n=\{\V^*\delta_t, t\in T, |t|\le n\}$. Then $|D_n|\le 2^{n+1}$ and
for any $t\in T$ we have
\[
dist(\V^*\delta_t, D_n)^2=
\begin{cases} 0, &|t|\le n, \\ 
             \sum_{l=n+1}^{|t|} w_l, & |t|>n.\end{cases}
\] 
We see that 
\[ 
   e_{n+2}(D)\le [ \sum_{l=n+1}^\infty w_l]^{1/2} = [\sum_{l=n+1}^\infty (1+l)^{-\beta}]^{1/2} 
   \le c \, n^{-\frac{\beta-1}2}.     
\]
Now recall that a polynomial upper bound $e_{n}(D)\le c n^{-\alpha}$ for any set $D$ in a Hilbert space 
yields a bound on $e_n(aco D)$, 
where $aco \, D$ denotes the absolutely convex hull of $D$. Namely, as established in
\cite{CKP} for $\alpha\not = 1/2$ and in \cite{Gao} for $\alpha=1/2$ under this assumption we have
\begin{eqnarray} \nonumber
e_n(aco \, D) &\le& C n^{-\alpha}, \qquad \qquad \qquad \quad \, \alpha<1/2,
\\ \label{enAacoA}
e_n(aco \, D) &\le& C n^{-1/2} \ln n, \qquad \quad \quad \ \alpha=1/2,
\\ \nonumber
e_n(aco \, D) &\le& C n^{-1/2} (\ln n)^{1/2-\alpha}\qquad \alpha>1/2.
\end{eqnarray}
By letting here $\alpha=\frac{\beta-1}2$, we obtain the desired upper bounds in Theorem \ref{t2noncrit},
because by the property of the unit ball in $\ell_1$-space we have $e_n(\V^*)=e_n(aco\, D)$.

\noindent {\it Lower bound}. For any $n\in \N$ let $m=2^n$ and denote 
$
  \{t:|t|=n \} := (t_j)_{1\le j\le m} \ .
$
Take any $(s_j)_{1\le j\le m}$ such that $|s_j|=2n$ and $s_j$ is an (indirect) offspring of $t_j$. 
Let $\mu_j=\delta_{s_j}-\delta_{t_j}$. Then
\[
(\V^*\mu_j))(u) = 1_{\{t_j\prec u \preceq s_j\}}.
\]
These image vectors are orthogonal, since they have disjoint supports, and for appropriate $C_1$
\[
|| \V^*\mu_j||_{2,W}^2  = \sum_{l=n+1}^{2n}w_l \ge C_1^2 \, n^{-(\beta-1)}.
\]
We notice that we found $m=2^n$ elements $\mu_j$ such that $||\mu_j||_1=2$
and for $i\not=j$ we have 
$ ||\V^*(\mu_j-\mu_i)||_{2,W}\ge C_1\, n^{-\frac{\beta-1}2}$. It follows
that
\[ 
  e_{n+1}(\V^*) \ge \frac{C_1}{2} \ n^{-\frac{\beta-1}2}.
\]
This is true for any $\beta>1$ but it is optimal only for $1<\beta\le 2$, while for $\beta>2$
we need a refined argument. 

By using the same vectors, we see that the restriction of $\V^*$ on the span of vectors
$(\mu_j)$ is isometric to the embedding $I_m:\ell_1^m\to\ell_2^m$ up to the coefficient 
\[
\frac{||\V^* \mu_j||_{2,W}}   {||\mu_j||_1} \ge
\frac{C_1}{2}\ n^{-\frac{\beta-1}2}. 
\]
Recall that with appropriate numerical $c>0$ we have a (sharp) estimate
\[
   e_k(I_m) \ge \left[\frac{c\ln(1+m/k)}{k}\right]^{1/2}, \qquad \log_2 m\le k \le m, 
\]
see \cite{Sc}.
Choose $n=n(k)$ such that $2^{n/2}\le k \le 2^{(n+1)/2}$.
Since $m=m(n)=2^n$, the parameter  $k$ fits in the range and we obtain
a bound that is sharp for $\beta\ge 2$,
\begin{eqnarray*}
  e_{k}(\V^*) &\ge& \frac{C_1}{2} \ n^{-\frac{\beta-1}2} e_k(I_m) 
  \\
  &\ge&  \frac{C_1}{2} \   \ n^{-\frac{\beta-1}2} \left[\frac{c\ln(1+2^{(n-1)/2})}{k}\right]^{1/2} 
  \ge  \frac{\tilde C_1(\ln k)^{1-\beta/2}}{k^{1/2}} \ .\ \Box
\end{eqnarray*}
\medskip

There are many available proofs for upper bounds in Theorem $\ref{t2noncrit}$. The one presented
here is probably the shortest one. It is due to W. Linde.
We refer to $\cite{AL}$ for the studies of other summation operators 
in probabilistic language.

\subsection{Critical case}
 
We see that in Theorem \ref{t2noncrit}  the upper and lower estimates for $\beta\not =2$ 
are of the same order and 
are "easy" to obtain modulo known results, although at the point $\beta=2$ the behavior
of entropy undergoes a striking change. Moreover, in the case $\beta=2$ the estimates of
Theorem \ref{t2noncrit} do not fit together and leave a logarithmic gap. 
Apparently this gap is impossible to close just by combining the existing results.
Therefore, we call the case $\beta=2$ a critical one. We will show that the lower bound of Theorem
\ref{t2noncrit} is in fact sharp but the proof of the corresponding upper bound
is by far more complicated and requires a new method.
 
\begin{thm} \label{t2} Let the weight $W$ be given by $(\ref{bw})$ with $\beta=2$. 
Consider the linear operator
$\V: \ell_2(T,W)\to \ell_\infty(T)$
defined by $(\ref{v_op})$
and its dual 
$\V^*:\ell_1(T)\to \ell_2(T,W)$
defined by $(\ref{vstar_op})$. There exists a numeric constant $C$ such that
 for all positive integers $n$ we have the following upper bound for its entropy
 numbers
\[
    \max \left \{e_n(\V),\, e_n(\V^*) \right\} \le \frac{C}{\sqrt{n}}\ . 
\]
\end{thm} 

{\bf Proof.} The proof of Theorem \ref{t2} will be splitted in few steps, each step having its clear own
meaning.
We keep the notation $\OO(t), \OO_l(t), |t|, s_\mu(t),||\mu||(t)$ from the previous section.

\subsubsection*{Step 1: essential subtrees}

A subset $\Upsilon\subset T$ is called a subtree, if $t\prec u \in \Upsilon$ yields $t\in \Upsilon$.
In particular, the root is contained in any subtree.

The evaluation of entropy numbers will be  based on the construction of some family of subtrees $\Upsilon^\mu$ 
based on a stopping rule.
Namely, let  $\sigma_l= \frac{l}{n}$.
For $\mu\in \ell_1(T)$ satisfying $||\mu||_1\le 1$
we define the {\it $n$-essential subtree} $\Upsilon^\mu$ by starting from the root and including
all nodes in $\Upsilon^\mu$ while $||\mu||(t) > \sigma_{|t|}$ and stopping the construction
while $||\mu||(t) \le \sigma_{|t|}$. 
We denote by $B^\mu$ the set of nodes where construction was stopped. 

Since $\sigma_{n+1}>1$, 
we have $\Upsilon^\mu\cap T_{n+1}=\emptyset$, that is we stop the 
construction not later than  at the level $n$. In particular, $\Upsilon^\mu$ is finite.
Notice that we have a partition
\begin{equation}\label{partition}
T=\Upsilon^\mu \bigcup \ \left( \bigcup_{t\in B^\mu} \OO(t) \right).
\end{equation}

Now we evaluate the size of $\Upsilon^\mu$ and will see that:
\smallskip

{\it The size of $\Upsilon^\mu$ is dramatically small.} 
\smallskip

This is a decisive step towards our goal.
Let $N_l=|\Upsilon^\mu\cap T_l|$. We have

\begin{lem} \label{lb1} Let $Q$ be the set of terminal nodes of $\Upsilon^\mu$.
It is true that
\begin{equation}\label{bsizeTmuQ}
\sum_{t\in Q} |t| \le n.
\end{equation}
Moreover,
\begin{equation}\label{bsizeTmu}
\sum_{l=1}^\infty N_l \le n +1.
\end{equation}
\end{lem}

\noindent {\bf Proof of Lemma \ref{lb1}.}\ 
By the definition of $n$-essential tree we have
\[
1 \ge \sum_{t\in Q} ||\mu||(t) \ge  \sum_{t\in Q} \frac{|t|}{n}\ .
\]
It follows that $\sum_{t\in Q} |t| \le n$, as required in (\ref{bsizeTmuQ}).
On the other hand, for any tree $\Upsilon$ and its terminal set $Q$ it is true that
\begin{eqnarray*}
\sum_l |\Upsilon_l|  &=& 1+ \sum_{u\in \Upsilon, |u|>0} 1  
\le 1+ \sum_{u\in \Upsilon, |u|>0} |Q\cap \OO(u)| 
\\
&=&
1 +  \sum_{t\in Q} \ \sum_{u: |u|>0,u \preceq t} 1 =
1 +  \sum_{t\in Q} |t|,
\end{eqnarray*}
thus (\ref{bsizeTmu}) follows. $\Box$ 
%%\medskip

%%\begin{rem} Actually, we proved something stronger than $(\ref{bsizeTmu})$, namely,
%%\[
%%\sum_{u\in \Upsilon^\mu, |u|>0} |Q\cap \OO(u)| \le n.
%%\]
%%\end{rem}

We can also easily evaluate the number of possible $n$-essential trees.

\begin{lem} \label{lb2} \ The number of subtrees of binary tree 
whose terminal set $Q$ satisfies $(\ref{bsizeTmuQ})$ does not exceed
$(4e)^n$.
\end{lem}

\noindent {\bf Proof of  Lemma \ref{lb2}.} \  
Since a subtree is entirely defined by its terminal set, we have to find out how many sets $Q$
satisfy (\ref{bsizeTmuQ}). Denote $q_l=|Q\cap T_l|$. Then (\ref{bsizeTmuQ}) writes as
\[
\sum_l l\, q_l \le n.
\]
Since $q_l\le \frac nl$, the number of non-negative integer solutions of this inequality does not exceed
\[
\prod_{l=1}^n (1+\frac nl) \le \prod_{l=1}^n \frac {2n}l =\frac{(2n)^n}{n!}
\le \frac{(2n)^n}{(n/e)^n}= (2e)^n.
\]
Moreover, for given sequence $q_l$, while constructing a set $Q$, on each level $l$ of a binary tree
we have to choose $q_l$ elements
from at most $2^{l}$ elements of this level. Therefore the number of possible sets not exceed
\[ 
\prod_{l=1}^n \left( { {2^{l}} \atop  {q_l} } \right) \le \prod_{l=1}^n (2^{l})^{q_l}
= 2^{\sum_{l=1}^n l q_l} \le 2^{n}. 
\]
 $\Box$

We finish the discussion of essential subtrees by proving their useful approximation property.
It follows from  (\ref{s1}) and (\ref{s2}) that for any $t\in T$ it is true that
\begin{eqnarray} \nonumber
   \sum_{l=|t|}^\infty \sum_{u\in \OO_l(t)} w(u) s_\mu(u)^2
   &\le& \sum_{l=|t|}^\infty \max_{u\in \OO_l(t)} |s_\mu(u)| \sum_{u\in \OO_l(t)} w(u) |s_\mu(u)|
   \\ \nonumber
   &\le& \sum_{l=|t|}^\infty \max_{u\in \OO_l(t)} ||\mu||(u) \sum_{u\in \OO_l(t)} w(u) ||\mu(u)||
  \\ \label{bs3}
   &\le&  ||\mu||(t)^2 \sum_{l=|t|}^\infty (1+|l|)^{-2}
   \le ||\mu||(t)^2  |t|^{-1}.
\end{eqnarray}

Moreover, (\ref{bs3}) and the definition of $B^\mu$ yield
\begin{eqnarray} \nonumber
&&\sum_{t\not\in \Upsilon^\mu}    s_\mu(u)^2 w(u) 
=
\sum_{t\in B^\mu}  \sum_{u\in \OO_l(t)}  s_\mu(u)^2 w(u)
\le
 \sum_{t\in B^\mu} ||\mu||(t)^2 |t|^{-1}
\\ \label{bs4}
&\le&   \sum_{t\in B^\mu} \sigma_{|t|} ||\mu||(t) |t|^{-1}
=  \sum_{t\in B^\mu} \frac{|t| }{ n}\ ||\mu||(t)\ |t|^{-1} 
=  \sum_{t\in B^\mu}  ||\mu||(t)\  n^{-1}
 \le   n^{-1}. 
\end{eqnarray}

Therefore, as we will see soon, the part of operator $\V^*$ 
related to the complement of $\Upsilon^\mu$ is not essential
at the precision level $n^{-1/2}$ which explains the name "essential" we gave to this family.

\subsubsection*{Step 2: approximating operators}

We are going now to construct a family of finite rank operators approximating the operator $\V^*$.
Each operator will correspond to an $n$-essential subtree. However the construction is valid for
any subtree of $T$. Given a subtree $\Upsilon\subset T$ we define three operators related to $\Upsilon$.
The operator $\V_\Upsilon^*:\ell_1(\Upsilon)\to \ell_2(\Upsilon,W)$ is
defined by
\[ (\V_\Upsilon^*\mu)(t)=\sum_{u\in \OO(t)\cap\Upsilon}  \mu(u), \qquad t\in \Upsilon.
\]
This is essentially the same operator as $\V^*$ restricted on elements supported by $\Upsilon$.

Now define a mapping $z$ from the complement of $\Upsilon$ to the boundary of $\Upsilon$
by letting $z(s)$ be the last node in $\Upsilon$ on the way from the root to $s$.
We denote $Z(t)=z^{-1}(t)$. This set will be non-empty only if $t$ belongs to the boundary
of $\Upsilon$. 

Now the flush-projection operator
$P_\Upsilon^*:\ell_1(T)\to \ell_1(\Upsilon)$ is
defined by
\[ (P_\Upsilon^*\mu)(t)= \mu(t)+ \sum_{u\in Z(t)}  \mu(u), \qquad t\in \Upsilon.
\]
This operator projects measures supported by $T$ onto the measures supported by $\Upsilon$.
It is clear that $||P_\Upsilon||\le 1$, i.e. $P_\Upsilon$ is a contraction.

The main property of the operators introduced so far reads as
\begin{equation} \label{flush}
    (\V_\Upsilon^* P_\Upsilon \mu)(t)=  \sum_{u\in \OO(t)}  \mu(u)=  (\V^* \mu)(t),
     \qquad \forall t\in \Upsilon, \mu\in \ell_1(T).
\end{equation} 

Finally, we will use the natural embedding $\iota_\Upsilon: \ell_2(\Upsilon,W) \to \ell_2(T,W)$
defined by
\[ 
\iota_\Upsilon \mu(t)=
\begin{cases}
\mu(t), & t\in \Upsilon, \\
0, & t\not \in \Upsilon.
\end{cases}
\]
Combining all together we define the approximating operator $\A_\Upsilon:\ell_1(T)\to \ell_2(T,W)$ by
$\A_\Upsilon=\iota_\Upsilon \V_\Upsilon^*P_\Upsilon$. It follows from (\ref{flush})
that for any subtree $\Upsilon$ and any $\mu\in \ell_1(T)$ we have
\[
    (\V^* \mu)(t) - (\A_\Upsilon) \mu(t) =  
    \begin{cases} 0, & t\in \Upsilon, \\
                  s_\mu(t), & t\not\in\Upsilon.
    \end{cases}              
\]
Hence,
\begin{equation} \label{approxA}
    ||(\V^* - \A_\Upsilon) \mu||^2_{2,W} =  
     \sum_{t\not\in\Upsilon}  |s_\mu(t)|^2 w_{|t|}.              
\end{equation} 

Finally notice that  since $||P_\Upsilon||\le 1$ and $||\iota_\Upsilon||\le 1$,
we have for any $m\in \N$
\begin{equation} \label{approxAem}
    e_m\left(\A_\Upsilon\right) \le e_m\left(\V_\Upsilon^*\right).              
\end{equation}

So far we have not specified our subtree. Now we will use the $n$-essential subtrees constructed above.
For any given $n$ let
 $\Gamma=\{\Upsilon\}$ be the set of all subtrees $\Upsilon\subset T$ satisfying
(\ref{bsizeTmuQ}), hence (\ref{bsizeTmu}). Recall that by Lemma \ref{lb2}  we have $|\Gamma| \le (4e)^n$ 
and for any $\mu\in \ell_1(T)$ its $n$-essential subtree $\Upsilon^\mu$ belongs to $\Gamma$.

By comparing inequality (\ref{bs4}) with (\ref{approxA}) we see that for any $\mu\in \ell_1(T)$ with
$||\mu||_1\le 1$ we have
\[
    ||(\V^* - \A_{\Upsilon^\mu}) \mu||^2_{2,W} \le n^{-1}.              
\]
In other words,
\begin{equation} \label{approxAmu}
    \sup_{\mu:||\mu||_1\le 1} \ \inf_{\Upsilon\in\Gamma} \ ||(\V^* - \A_\Upsilon) \mu||_{2,W} \le n^{-1/2}.              
\end{equation}
Recall that for every $\mu$ its own approximating operator is used. We will show now how
the properties like this one can be applied. This simple idea seems to be of independent interest, thus
we state it as a separate statement.

\subsubsection*{Step 3: approximation lemma}

The following lemma shows how a linear operator $V$ can be approximated by a family
of operators $(V_\gamma)_{\gamma\in\Gamma}$ in a sense that for every element $x$ its image $Vx$ is approximated
by $V_\gamma x$ with appropriate $\gamma$ depending of $x$.  

\begin{lem} 
\label{approxV}
Let $X,Y$ be the normed spaces and $V,\ (V_\gamma)_{\gamma\in \Gamma}$ be the
linear operators acting from $X$ to $Y$. Then for any $n\in \N$ it is true that
\begin{equation}
\label{s1s2} 
e_{n+[\log_2 |\Gamma|]+1} (V) \le  
\sup_{\gamma\in \Gamma} e_n(V_\gamma) + \sup_{x\in B_X} \inf_{\gamma\in \Gamma} ||Vx-V_\gamma x||_Y ,
\end{equation}
where $B_X=\{x\in X:\ ||x||_X\le 1\}$.
\end{lem}

{\bf Proof:\ }
Denote $S_1$ and $S_2$ the expressions in (\ref{s1s2})
and fix a small $\delta>0$.
For every $\gamma$ we can choose an $(e_n(V_\gamma)+\delta)$-net $N_\gamma$ of size $2^{n-1}$ for the set
$V_\gamma(B_X)$ in the space $Y$. Let
\[ 
N= \bigcup_{\gamma\in \Gamma} N_\gamma
\]
be a global net. Clearly, 
\[ 
  \#\{N \} = |\Gamma| \cdot 2^{n-1} \le 2^{[\log_2 |\Gamma|]+n}  .
\]
For any $x\in B_X$ we first find a $\gamma$ such that
\[
   ||Vx-V_\gamma x||_Y \le S_2 + \delta.
\]
Then we find an element $y\in N_\gamma\subset N$ such that
\[
   ||V_\gamma x-y||_Y \le e_n(V_\gamma)+\delta.
\]
By triangle inequality, we have
\[
   ||Vx- y||_Y \le S_2 + \delta +e_n(V_\gamma)+\delta\le S_1+S_2 +2\delta.
\]
Therefore, $N$ is an $(S_1+S_2 +2\delta)$-net for the set $V(B_X)$ and its size does
not exceed $2^{[\log_2 |\Gamma|]+n}$.
The assertion of lemma follows.
$\Box$

We will apply Lemma \ref{approxV} to $X=\ell_1(T)$, $B_X=\{\mu:||\mu||_1\le 1\}$, 
$Y=\ell_2(T,W)$, $V=\V^*$, and approximating family of operators $(\A_\Upsilon)_{\Upsilon\in \Gamma}$. Now
(\ref{s1s2}), (\ref{approxAmu}), and (\ref{approxAem})  together with the known estimate of $|\Gamma|$ 
given in Lemma \ref{lb2} yield

\begin{equation}
\label{log4e} 
e_{(1+[\log_2(4e) ])n+1} (\V^*) \ \le\   
\sup_{\Upsilon\in \Gamma} e_n(\V_\Upsilon^*) +  n^{-1/2}. 
\end{equation} 

\subsubsection*{Step 4: evaluation of operators on short trees}

With (\ref{log4e}) at hand, it remains to evaluate $e_n(\V_\Upsilon^*)$ 
for fixed $\Upsilon\in \Gamma$. In other words,
we have to evaluate the entropy of the operator restricted 
to a tree of a very small size (due the bound 
(\ref{bsizeTmu}) for the size of $\Upsilon$).

Towards this aim, recall an important entropy bound from \cite{CP}, Corollary 2.4\, (i). 
There exists a constant $c>0$ such that for any operator ${\mathcal W}$ acting from
$\ell_1^m$ to a Hilbert space and any $k\in \N$ it is true that
\begin{equation}
\label{Carl1}
     e_k({\mathcal W}) \le c \, \ln^{1/2}(m+1) \, ||{\mathcal W}|| \, k^{-1/2} .
\end{equation}

We will apply this estimate to particular situation of tree operators. 

Let $\Delta$ be the tree that consists of the first 
$[\ln n/4]$ levels of binary tree.

Let us split our operator in a sum $\V_\Upsilon^*=\V_\Upsilon^+ + \V_\Upsilon^0$, where
$\V_\Upsilon^+$ corresponds to the layers distant from the root,
\[
\V_\Upsilon^+ \mu (t)=
\begin{cases} 0,& t\in\Delta,
\\
\V_\Upsilon^*\mu(t),& t\not \in\Delta,
\end{cases}
\] 
while $\V_\Upsilon^0$ corresponds to the first layers
\[
\V_\Upsilon^0 \mu (t)=
\begin{cases} \V_\Upsilon^*\mu(t),& t\in\Delta,
\\
0,& t\not \in\Delta.
\end{cases}
\] 

%%\[
%%w_l=
%%\begin{cases} (1+l)^{-2},& l\ge \ln n/4,
%%\\
%%0,& l<\ln n/4.
%%\end{cases}
%%\] 
%%while $\V_\Upsilon^0$ corresponds to
%%\[
%%w_l=
%%\begin{cases} 0,& l\ge \ln n/4,
%%\\
%%(1+l)^{-2},& l<\ln n/4.
%%\end{cases}
%%\]

The idea behind this splitting is simple: the operator $\V_\Upsilon^+$ has a 
small norm while $\V_\Upsilon^{0}$ has a small image dimension.
We first study the operator $\V_\Upsilon^+$.
Notice that 
\[
   ||\V_\Upsilon^+|| \le \left(\sum_{l\ge \ln n/4} (1+l)^{-2}\right)^{1/2} 
   \le (\ln n/4)^{-1/2}= 2(\ln n)^{-1/2}.
\]
For any tree $\Upsilon$ of size bounded by $m$,  from (\ref{Carl1}) 
we get 
\begin{equation}
\label{Carl2}
           e_k(\V_\Upsilon^+) \le c \, \ln^{1/2}(m+1) \, ||\V_\Upsilon^+|| \, k^{-1/2}. 
\end{equation}
and applying this  with $k=n$, $m=n+1$  we obtain
\[
      e_n(\V_\Upsilon^+) \le c  \, n^{-1/2}. 
\]

Now we have to consider the operator $\V_\Upsilon^0$.

Notice that since weights on higher
levels vanish, operator $\V_\Upsilon^0$ actually acts into $\ell_{2, W}(\Delta)$. 
The size of $\Delta$ is merely $2^{1+\ln n/4}\le 2n^{1/4}$,
thus estimation can be rather crude.

Write $\V_\Upsilon^0 = I \circ \V_\Upsilon^{00}$, where $\V_\Upsilon^{00}$ is the same operator
as $\V_\Upsilon^0$ but  acting into $\ell_\infty(\Delta)$ and $I$ is the embedding
of $\ell_\infty(\Delta)$ in $\ell_{2, W}(\Delta)$. The operator $\V_\Upsilon^{00}$
is a contraction, since
\[
   ||\V_\Upsilon^{00}\mu||_\infty
= \max_{t\in\Delta} \left| \sum_{u\in\OO(t)\cap\Upsilon}   \mu(u) \right| \le ||\mu||_1.
\]

On the other hand,  we can easily evaluate the entropy of $I$.
The net $\HH_\Delta\subset \ell_{2, W}(\Delta)$ 
will consist of all possible functions $h$ of the form
\[
h(t)= j(t) n^{-1}, \qquad t\in \Delta,
\]
where $j(t)$ are odd integers satisfying $|j(t)| \le n$. 
Notice that there are no more than $2n$
choices for each $j(t)$. 

Now we provide the estimates for approximation error and for the size of $\HH_\Delta$.
We start with evaluating approximation error. Let $x\in \ell_\infty(\Delta)$ be such that
$||x||_\infty\le 1$. Then for any $t\in \Delta$ we have $|x(t)|\le 1$, hence, there exists a function
$h\in \HH_\Delta$ such that
\[
|x(t)-h(t)|\le n^{-1}, \qquad \forall t\in \Delta.
\]
Therefore,
\[
||Ix-h||_{2,W}^2 = \sum_{t\in \Delta} w_{|t|} |x(t)-h(t)|^2
\le  |\Delta| n^{-2} \le 2n^{-7/4}.
\]
The size of $\HH_\Delta$ is bounded by
\[
|\HH_\Delta|\le (2n)^{|\Delta|} \le (2n)^{2n^{1/4}}
=2^{2n^{1/4}(1+\log_2 n)}\le 2^{2(n+1)}.
\]
We conclude that
\[
e_{2n+3}(\V_\Upsilon^0)\le e_{2n+3}(I)\le 2n^{-7/8},
\]
and we are done with operator $\V_\Upsilon^0$, too. 

Having the bounds both for and $e_{n}(\V_\Upsilon^+)$ and $e_{n}(\V_\Upsilon^0)$
by standard entropy estimates we get a bound for the sum of operators, i.e.  
\[
   e_{n}(\V_\Upsilon^*) \le c\ n^{-1/2}.
\]
Finally, it follows from (\ref{log4e}) that
\[
   e_{n}(\V^*) \le c\ n^{-1/2},
\]
as required in the assertion of Theorem \ref{t2}. 
Once the bound for $e_n(\V^*) $ is obtained,
 the bound for $e_n(\V)$ follows from famous duality connection
for entropy numbers,
\begin{equation}\label{duality}
    e_n(\V) \le c_1 \ e_{c_2n}(\V^*) 
\end{equation}
for some numerical constants $c_1,c_2$,  which is still a conjecture for general Banach spaces but is a proved 
statement in our situation (one of the spaces is a Hilbert one), see 
%%\cite{AMS1}, 
\cite{AMS2}.
%%Actually, an older result \cite{TJ} would suffice.\
$\Box$

%%%%%%%%%%%%%%%%%%%%%%%%%%%%%%%%%%%%%%%%%%%%%%%%%%%%%%%%%%%%%%%%%%%%%%%%%%%%%%%%%%%%%%%

\section{Entropy of an integral operator }
\label{s:integral}
Let $r<e^{-2}$ be a small number.
In this section $(\cdot,\cdot)$ and $||\cdot||$ denote the scalar product 
and the norm in $L_2[0,r]$, respectively.
We denote by $\MM[0,r]$ the space of signed measures of finite variation and
$||\cdot||_1$ the respective variation norm. Moreover, $||\mu||_1(I)$
stands for the variation of $\mu \in \MM[0,r]$ on an interval $I$.

Our aim is to study the critical integral operator 
$\V:L_2[0,r]\to\CC[0,r]$ defined by
\[
\V f(t)=\int_0^t f(s)K_t(s) ds =(f,K_t), \qquad 0\le t \le r,
\]
and its adjoint
$\V^*:\MM[0,r]\to L_2[0,r]$ defined by
\[
\V^* \mu(s)=\int_0^r K_t(s) \mu(dt), \qquad 0\le s \le r,
\]
where the critical kernel is
\[
 K_t(s)=(t-s)_+^{-1/2}|\ln(t-s)_+|^{-1}.
\]

Before we start the studies of $K$, let us explain why it is critical in our context. 
Consider the family of kernels
\[
 K^{(\beta)}_t(s)=(t-s)_+^{-1/2}|\ln(t-s)_+|^{-\beta}, \qquad 1/2<\beta <\infty,
\]
and the corresponding operators $\V_\beta$. It is known from the works of Linde and Lacey 
\cite{La}, \cite{Lin}  that
\begin{eqnarray*}
c\  n^{1/2-\beta} \ \le \ e_n(\V_\beta)\ &\le&  C\   n^{1/2-\beta}, \qquad  \quad \quad 1/2< \beta < 1,  
 \\ 
c\ n^{-1/2} \ \le \ e_n(\V_\beta)\   &\le&  C\ n^{-1/2} \ln n, \qquad \quad \, \beta = 1, 
\\ 
c\ n^{-1/2} (\ln n)^{1-\beta} \ \le \ e_n(\V_\beta)\   &\le&  C \  n^{-1/2} (\ln n)^{1-\beta}, \quad  \beta>1.
\end{eqnarray*}

Therefore, we see that the most interesting kernel $K=K^{(1)}$ lays on the boundary between two different regimes
and observe a logarithmic gap between the lower and upper bounds.
The situation is exactly the same as in Theorem \ref{t2noncrit}.

The main property of the kernel $K$ we need is its  modulus 
of continuity.\footnote{However, we will also use that 
the kernel $u^{-1/2}|\ln(u)|^{-1}$ is a decreasing convex
function on $[0,r]$ by the choice of $r$.}
An elementary calculation shows that
%%\begin{lem} \label{li1} 
for all $0\le t\le t+u\le r$ 
%% it is true that
\begin{equation} \label{li1}
   ||K_{t+u}-K_t||_2 \le 2 |\ln u|^{-1/2}.
\end{equation}
%%\end{lem}
%%
%%{\bf Proof of Lemma \ref{li1}.}\ We have
%%\begin{eqnarray*}
%%&& ||K_{t+u}-K_t||_2^2 
%%\\
%%&=& \int_0^t \left( (t-s)^{-1/2} |\ln(t-s)|^{-1} - (t+u-s)^{-1/2} |\ln(t+u-s)|^{-1}\right)^2ds
%%\\
%%&& + \int_t^{t+u} (t+u-s)^{-1} |\ln (t+u-s)|^{-1} ds 
%%\\
%%&=& \int_0^t \left( v^{-1/2} |\ln v|^{-1} - (v+u)^{-1/2} |\ln(v+u)|^{-1}\right)^2dv
%% + \int_0^{u} v^{-1} |\ln v|^{-2} dv 
%%\\
%%&\le& \int_u^t \left( v^{-1/2} |\ln v|^{-1} - (v+u)^{-1/2} |\ln(v+u)|^{-1}\right)^2dv
%% + 2 \int_0^{u} v^{-1} |\ln v|^{-2} dv 
%%\\
%%\\
%%&\le& u^2 \int_u^t  v^{-3} |\ln v|^{-2} dv
%% + 2 \int_0^{u} v^{-1} |\ln v|^{-2} dv 
%%\le  |\ln u|^{-2} + 2 |\ln u|^{-1}, 
%%\end{eqnarray*}
%%and the assertion of lemma follows. $\Box$

\begin{thm} \label{ti} 
For all positive integers $n$ and for a numerical constant $C$ we have
\[
\max \left\{ e_n(\V),  e_n(\V^*) \right\} \le   {C}   {n^{-1/2}}\ .
\]
\end{thm}

{\bf Proof of Theorem  \ref{ti}.} 
We repeat the ideas applied earlier to the summation operator on a binary tree.
We first find a family of good finite rank approximations to $\V^*$ by giving interpretation
for $n$-essential subtrees. We will construct {\it $n$-essential partition $\II^\mu_n$
of $ [0,r]$} as follows.  Given a positive integer $n$ and an element 
$\mu\in\MM[0,r]$ we start dividing the interval $[0,r]$ in halves and continue dividing
while a (binary) interval $I= \left(\frac{ir}{2^l},\frac{(i+1)r}{2^l}\right] $ subject to division satisfies 
\begin{equation}\label{critdiv}
                                    ||\mu||_1(I)\ge \frac l n.
\end{equation}
Once an interval does not satisfy (\ref{critdiv}) we do not divide it and include it in 
our partition $\II^\mu_n$. If $||\mu||_1\le 1$, the condition (\ref{critdiv}) fails for $l>n$.
Therefore, our construction provides a finite partition of $[0,r]$ in binary intervals
of variable length.

The partition $\II^\mu_n$ depends on $\mu$ but we will show now that the
number of possible partitions and their size are rather limited.

Let $\D$ be the set of all binary intervals we divided during the
construction of $\II^\mu_n$. Notice that $\D$ is a tree w.r.t. inclusion.
Let $Q$ be the set of all terminal intervals of $\D$. In other words,
$I\in Q$ iff $I$ satisfies (\ref{critdiv}) but neither of its halves 
satisfies it.  It is important for us that $Q$ uniquely determines both
$\D$ and $\II^\mu_n$. Indeed, any subtree of the binary tree is determined by the
set of its terminal nodes. Thus $Q$ determines $\D$. Moreover, $\II^\mu_n$ consists
of all direct offsprings of elements of $\D$ that do not belong to $\D$.

Let $q_l=\#\{ I\in Q:|I|=2^{-l}r\}$. Then by (\ref{critdiv})
\[
   1\ge ||\mu||_1\ge \sum_{I\in Q} ||\mu||_1(I)=
   \sum_{l=0}^\infty \sum_{I\in Q, |I|=2^{-l}} ||\mu||_1(I) 
\ge \sum_{l=0}^{\infty} q_l \frac ln\ .
\]
Hence,
\[
\sum_{l=0}^\infty q_l \cdot l \le n.
\]
By Lemma \ref{lb2}, the number of possible trees $Q$, thus the number of possible $n$-essential 
partitions does not exceed $(4e)^n$. It it is also worthwhile to notice that the
number of intervals in $\II^\mu_n$ satisfies 
\begin{equation}\label{sizeInmu}
|\II^\mu_n|\le 2|\D|\le 2(n+1)
\end{equation}
 by Lemma \ref{lb1}.

Consider a finite dimensional approximation for $\V^*$ generated by any partition $\II$, the operator
$\V^*_\II: \MM[0,r]\rightarrow L_2[0,r]$ defined by
\[
(\V^*_\II \mu)=\sum_{I\in\II} \mu(I) K_{t_I},
\]
where $t_I$ is the left end of $I$. We evaluate the approximation error 
$\Delta_\II=\V^*- \V^*_\II$. By the definition,
\[
(\Delta_\II\, \mu)=\sum_{I\in\II} \int_I \left(K_t-K_{t_I}\right) \mu(dt),
\]

We are going to show that the approximation error is particularly small when we use 
the $n$-essential partition.

\begin{prop} \label{p1i}
For any  $n\in \N$ and any $\mu$ with $||\mu||_1\le 1$ we have
\begin{equation}\label{mainconj}
    ||\Delta_{\II^\mu_n} \mu ||_2 \le C \ n^{-1/2}\ .
\end{equation}
\end{prop}

{\bf Proof of Proposition \ref{p1i}.}
Let $\mu=\mu_+ -  \mu_-$ be the Hahn decomposition of $\mu$. It is enough  to show that
\begin{equation}\label{mainconjp}
    ||\Delta_{\II^\mu_n} \mu_+ ||_2 \le  C \ n^{-1/2}\ 
\end{equation}
and to prove the similar inequality for $\mu_-$. We start with
\begin{eqnarray*}
     ||\Delta_{\II^\mu_n} \mu_+ ||_2  &=& \left(\Delta_{\II^\mu_n} \mu_+, \Delta_{\II^\mu_n} \mu_+\right)
 \\
    &=&
    \sum_{I_1, I_2\in \II^\mu_n }  
    \left(\int_{I_1} \left(K_t-K_{t_{I_1}}\right) \mu_+(dt),
          \int_{I_2} \left(K_t-K_{t_{I_2}}\right) \mu_+(dt)  \right)
\\
    &=&
   \sum_{I_1, I_2\in \II^\mu_n}  
   \int_{I_1} \int_{I_2} \left(K_{t_1}-K_{t_{I_1}}, K_{t_2}-K_{t_{I_2}} \right) 
   \mu_+(dt_1)\mu_+(dt_2). 
\end{eqnarray*}

For the main (diagonal) terms of this sum we have
\begin{eqnarray*} 
  &&  \sum_{I\in \II^\mu_n} 
  \int_I \int_I \left(K_{t_1}-K_{t_I}, K_{t_2}-K_{t_I} \right) \mu_+(dt_1)\mu_+(dt_2) 
   \\
  &\le& 
  \sum_{I\in \II^\mu_n} 
  \int_I \int_I ||K_{t_1}-K_{t_I}||_2 \, ||K_{t_2}-K_{t_I}||_2\ \mu_+(dt_1)\mu_+(dt_2)
   \\
  &\le&  
  \sum_{I\in \II^\mu_n} 
  \max_{t\in I} ||K_{t}-K_{t_I}||^2_2 \ \mu_+(I)^2
   \\
   &\le&  
  \sum_{I\in \II^\mu_n} 
   4 (\ln|I|)^{-1} \ \mu_+(I)^2 \qquad \textrm{by  (\ref{li1})}
   \\
  &\le& 
  \sum_l \sum_{I\in \II^\mu_n, |I|=2^{-l}r} 
   4 (\ln 2\cdot l)^{-1}\cdot  \frac{l} n \ \mu_+(I)\qquad \textrm{by definition of}\ \II^\mu_n
   \\
   &=&  \frac{4}{(\ln 2) \, n} \sum_{I\in \II^\mu_n} \mu_+(I) \le \frac{4}{(\ln 2)\, n}\ .
\end{eqnarray*}  
Unlike to the tree case, the summands in the definition of 
$\Delta_\II$ are not orthogonal, therefore
we can not stop here. We will show that the non-diagonal terms do not give a positive contribution
to the quantity we evaluate.
%% From the probabilistic viewpoint, this is very much in the spirit 
%% of negative dependence of the increments
%%for fractional Brownian motion with parameter $H<1/2$. We will state the corresponding result as an independent lemma.

Let $g:\rr\to \rr$ be a function such that $g$ vanishes on $(-\infty,0]$ and $g$ is a decreasing convex
non-negative function on $(0,+\infty)$. Let $K_t(\cdot)=g(t-\cdot)$ for $t\ge 0$.

\begin{lem} \label{li_negdep} For all $0\le a\le b\le c\le d\le r$ we have
\[
\int_0^r (K_d-K_c)(K_b-K_a) \le 0.
\]
\end{lem}

{\bf Proof of Lemma \ref{li_negdep}.} First of all, let us notice that the function 
$s\to K_d(s)-K_c(s)= g(d-s)-g(c-s)$ is non-positive and non-increasing while $s\in[0,c]\supset
[0,b]$.
Next, the function $s\to K_b(s)-K_a(s)= g(b-s)-g(a-s)$ is positive on $[a,b]$ and negative on $[0,a]$.
Therefore,
\begin{eqnarray*}
  \int_0^r (K_d-K_c)(K_b-K_a) &=& \int_0^b (K_d-K_c)(K_b-K_a) 
  \\
 &=&  \left( \int_0^a  + \int_a^b \right) (K_d-K_c)(K_b-K_a)  
 \\
 &\le&  (K_d(a)-K_c(a)) \left( \int_0^a  (K_b-K_a) + \int_a^b K_b \right) 
 \\
  &=&  (K_d(a)-K_c(a)) \left( \int_{b-a}^{b} g +  \int_{0}^{a} g  -  \int_{0}^{b-a} g \right) 
 \\
  &=&  (K_d(a)-K_c(a)) \int_{a}^{b} g \le 0. \ \Box
\end{eqnarray*}  
\medskip

By applying this result to our function $g(t)=t_+^{-1/2}|\ln t|^{-1}$ we obtain for any
$t_1\in I_1\in\II^\mu_n$,\ $t_2\in I_2\in\II^\mu_n$
\[
  \left(K_{t_1}-K_{t_{I_1}}, K_{t_2}-K_{t_{I_2}}\right) 
  =    \int_0^r (K_{t_1}-K_{t_{I_1}}) (K_{t_2}-K_{t_{I_2}}) \le 0,
\]
provided $I_1\not = I_2$.
Hence,
\[
\sum_{I_1, I_2\in \II^\mu_n, I_1\not= I_2} 
  \int_{I_1} \int_{I_2} \left(K_{t_1}-K_{t_{I_1}}, K_{t_2}-K_{t_{I_2}} \right) 
  \mu_+(dt_1)\mu_+(dt_2) \le 0. 
\]
Therefore,
\[
  ||\Delta_{\II^\mu_n} \mu_+ ||_2^2 
   \le 
  \sum_{I \in \II^\mu_n} 
  \int_{I} \int_{I} \left(K_{t_1}-K_{t_{I}}, K_{t_2}-K_{t_{I}} \right) 
  \mu_+(dt_1)\mu_+(dt_2)
   \le \frac{4}{(\ln 2)\, n}\ , 
\]  
and (\ref{mainconjp}) follows. The same inequality for $\mu_-$ is obtained by applying
(\ref{mainconjp}) to $-\mu$. Now (\ref{mainconj}) is proved completely. $\Box$
\medskip

We continue the proof of Theorem \ref{ti}. Let $\III_n=\{\II^\mu_n:||\mu||_1\le 1\}$
 be the set of all possible $n$-essential partitions of $[0,r]$.
 Recall that
 \begin{equation} \label{IIIsize}
 |\III_n|\le (4e)^n.
 \end{equation}
We claim that
\begin{equation} \label{Jn}
     \sup_{\II\in\III_n} e_n(\V^*_\II)\le {C}{n^{-1/2}}\ .  
\end{equation}  
Assuming this is obtained, the application of Lemma \ref{approxV} to the family of operators $\{\V^*_\II, \II\in \III_n\}$
along with the estimate of approximation error (\ref{mainconj}) and the estimate for the number
of operators (\ref{IIIsize})  lead to
$e_n(\V^*) \le {\tilde C}{n^{-1/2}}$
as required by assertion of Theorem \ref{ti}. 
The same estimate for $e_n(\V)$ follows by the duality argument (\ref{duality}).

Now it only remains to prove (\ref{Jn}).  Let us fix a partition $\II\in\III_n$. 
From now on, we do not need any particular properties of $n$-essential partitions, except
for the size bound (\ref{sizeInmu}). 

Consider an
auxiliary partition $\EE$ of $[0,r]$ constructed as follows. Take $m$ such that 
$2^{-m}\le n^{-1/4}\le 2^{1-m}$. Divide $[0,r]$ in binary intervals of length $r 2^{-m}$.
If a union of such intervals belongs to $\II$, then replace them by this union.
The result is a partition $\EE$. Notice that $\II$ is a refinement of $\EE$ and 
$|\EE|\le 2^m\le 2n^{1/4}$.
Write
\[
\V^*_\II=\V^*_\EE+(\V^*_\II-\V^*_\EE)
\]
and evaluate the entropy of both operators.

First we handle the low rank operator $\V^*_\EE$. 
%%Recall that
%%\[
%%\V^*_\EE \mu =\sum_{I\in\EE} K_{t_I} \mu(I).
%%\]
Consider the net
\[
\NN_\EE= \left\{ \sum_{I\in\EE} K_{t_I} {j_I}n^{-1},\ 
j_I\in \{1-n,\dots,0,\dots, n-1\}\right\}.
\]
Notice that
\[
|\NN_\EE|\le (2n-1)^{|\EE|}\le (2n)^{2n^{1/4}} \le 2^{2(n+1)}.
\]
On the other hand, for any $\mu$ with $||\mu||_1\le 1$ find an 
$h= \sum_{I\in\EE} K_{t_I} \frac{j_I}n$ such that
$\max_{I\in\EE} |\mu(I)-\frac{j_I}n| \le n^{-1}$. We have
\[
||\V^*_\EE \mu -h||_2 \le \sum_{I\in\EE} ||K_{t_I}||_2 \cdot |\frac{j_I}n-\mu(I)|
 \le \max_{t\in[0,1]} ||K_{t}||_2 \cdot n^{-1}\cdot |\EE|\le C \, n^{-3/4}.
\]
It follows that
\[
e_{2n+3}(\V^*_\EE )\le C \, n^{-3/4}.
\]
Now we handle the operator $\V^*_\II-\V^*_\EE$ which has a larger rank but smaller norm. 
By the definition,
\begin{eqnarray*}
  (\V^*_\II-\V^*_\EE)\mu &=&  \sum_{J\in\II} K_{t_J} \mu(J) - \sum_{I\in\EE} K_{t_I} \mu(I)
\\
   &=&  \sum_{J\in\II \backslash \EE} K_{t_J} \mu(J) - 
    \sum_{I\in\EE \backslash \II} K_{t_I} \mu(I)
\\
  &=&  \sum_{I\in\EE \backslash \II} \sum_{J\in\II,J\subset I } (K_{t_J}-K_{t_I})\mu(J).
\end{eqnarray*}  
Notice that the conditions $I\in\EE \backslash \II, J\subset I$ imply
\[
|t_J-t_I|\le |I|=2^{-m}\le n^{-1/4},
\]
hence by Lemma \ref{li1}
\[
||(\V^*_\II-\V^*_\EE)\mu||_2 \le 2 (\ln (n^{1/4}))^{-1/2} \sum_{J\in \II} |\mu(J)|
\le (\ln n)^{-1/2} ||\mu||_1,
\]
which simply means that  $||\V^*_\II-\V^*_\EE||\le (\ln n)^{-1/2}$.

On the other hand, recall that by (\ref{sizeInmu})
\[
\rank(\V^*_\II-\V^*_\EE)\le |\II|\le 2(n+1).
\]

We apply, as we did in the investigation of tree summation, the estimate (\ref{Carl1}) 
and obtain
\begin{eqnarray*}
e_k(\V^*_\II-\V^*_\EE) &\le& c \ \ln^{1/2} (\rank(\V^*_\II-\V^*_\EE)+1)\  ||\V^*_\II-\V^*_\EE||\ k^{-1/2}
\\ 
&\le& c \ \ln^{1/2}(2n+3) \cdot  (\ln n)^{-1/2} \cdot k^{-1/2} 
\\
&\le& c \ k^{-1/2}.
\end{eqnarray*}
By letting $k=n$,
\[
e_n(\V^*_\II-\V^*_\EE) \le  c \ n^{-1/2}.
\]
We conclude that
\[
e_{3n+1}(\V^*_\II) \le e_{2n+2}(\V^*_\EE)+ e_n\left(\V^*_\II-\V^*_\EE\right)
\le c\, n^{-1/2},
\]
and  (\ref{Jn}) follows. $\Box$
\medskip

%\section{%
{\bf Relation to the entropy of convex hulls}.
Recall a well known problem from the geometry of Banach spaces, see e.g. \cite{CKP}, 
\cite{LiLin}. 
Let $A$ be a set in a Hilbert space and
$aco\,A$ its absolutely convex hull. If we know the behavior of entropy numbers $e_n(A)$,
what can we say about $e_n(aco\,A)$? In (\ref{enAacoA})
we already recalled some known relations.
%%, it is known
%%that for $\beta<1/2$
%%\[  e_n(A)\le c \, n^{-\beta} \ \textrm{implies}\  e_n(aco\,A)\le C\, n^{-\beta},
%%\] 
%%while
%%\[  e_n(A)\le c\, n^{-1/2} \ \textrm{only implies}\  e_n(aco\,A)\le C\, n^{-1/2}\ln n.
%%\]
F.Gao \cite{Gao} was the first to construct a critical set $A$ with properties
\begin{equation}\label{Gaoset}  
   e_n(A)\le c\, n^{-1/2} \ \textrm{and} \   e_n(aco\,A)\ge C\, n^{-1/2}\ln n.
\end{equation}
We can call {\it Gao set} any set satisfying (\ref{Gaoset}). Later on, his arguments were streamlined and extended
to non-Hilbert case 
%%by J. Creutzig and I. Steinwart 
in \cite{CrSt}.  

The relation to our problem is the following. 
%%For simplicity, let us take a setting of summation
%%operator on the binary tree (Section \ref{s:binary}). 
Consider the critical tree summation operator $\V^*$ with the weight (\ref{bw}) and $\beta=2$. Take a set
$
   A=\left \{ \V^*(1_{\{t\}}), \quad t\in T \right\}
$ 
in the Hilbert space $\ell_2(T,W)$. It is plain that $e_n(A) \le C n^{-1/2}$, 
hence
%%and we see from the general fact above
that $e_n(aco\,A)\le C \, n^{-1/2}\ln n$. Since it is quite difficult to get a better upper bound,
one could think of $A$ as a candidate to be a Gao set, although of a nature very different from 
the known ones.
On the other hand,  $aco\,A$ is the image of the unit ball w.r.t. operator $\V^*$. In other words,
$e_n(\V^*)=e_n(aco\,A)$. Therefore, Theorem \ref{t2} shows that $A$ {\it is not a Gao set}.

%%Curiously enough, one can construct another set $B=\{b_t,t\in T\}$, such that $B$ is a Gao set, while
%%the distance between any pair of points of $B$ is smaller than that between their counterparts in $A$.
%%There is no contradiction, of course, but the reader can see from this observation
%%that the problem is really delicate.
%%\medskip

{\bf Acknowledgements}. I am grateful to M. Lacey and W. Linde for the problem statement, 
for interesting discussions around it and for informing me about their unpublished results.  
 Sincere thanks are also due to American Institute of Mathematics %%(Palo-Alto) and its staff  
 for organizing
 the related conference "Small ball inequalities in analysis, probability, and irregularities of distribution".
 %%(December 2008) where many people mentioned here had a chance to discuss
 %%this and related problems.  


\begin{thebibliography}{1000}
{\baselineskip=12pt \parskip=-2pt

%%\bibitem{AMS1} Artstein S., Milman V.D., Szarek S.J. (2003) Duality of metric entropy
%%in Euclidean case. C. R. Acad. Sci. Paris, S\'er. I, 337, 711--744.

\bibitem{AMS2} Artstein S., Milman V.D., Szarek S.J. (2004) Duality of metric entropy. 
Annals of Math., 159, 1313--1328.
\bibitem{AL} Aurzada F., Lifshits M. (2008) Small deviation probability via chaining,  
Stoch. Proc. Appl., 118, 2344--2368.
\bibitem{CKP}
Carl B., Kyrezi I., Pajor A. (2000) Metric entropy of convex hulls
in Banach spaces. J. London Math. Soc., 60, 871--896.
\bibitem{CP}
Carl B., Pajor A. (1988)
Gelfand numbers of operators with values in a Hilbert space.
Invent. Math., 94, 479--504.
\bibitem{CS} Carl B., Stephani I. (1990)
Entropy, Compactness and Approximation of Operators.
Cambridge University Press.
\bibitem{CrSt} Creutzig J., Steinwart, I. (2002) Metric entropy of convex hulls in type $p$ spaces --
the critical case. Proc. Amer. Math. Soc., 130, 3, 733--743.
\bibitem{ET} Edmunds D.E., Triebel H. (1996) Function Spaces, Entropy
Numbers and Differential Operators. Cambridge University Press.
\bibitem{Gao} Gao F. (2001) Metric entropy of convex hulls.
Isr. J. Math., 123, 359--364.
\bibitem{La} Lacey M. (2008) Private communication. 
\bibitem{LiLin} Li W.V., Linde W. (2000) Metric entropy of convex hulls in Hilbert spaces. Studia Math., 139, 29--45.
\bibitem{Lin} Linde W. (2008) Nondeterminism of linear operators and lower entropy estimates. J. Fourier Anal. Appl.,
14, 568--587.
\bibitem{Sc} Sch\"utt C. (1984) Entropy numbers of diagonal operators
between symmetric Banach spaces. J. Approx. Theory, 40, 2, 121-128.
%%\bibitem{TJ} Tomczak-Jaegermann N. (1987) Dualit\'e des nombres d'entropie pour les op\'erateurs \`a valeurs 
%%dans un espace de Hilbert. C. R. Acad. Sci. Paris, S\'er. I, 305, 299-301.
}
\end{thebibliography}
\end{document}